\documentclass[12pt,reqno]{amsart}
\usepackage{amssymb}
\usepackage[a4paper,margin=3cm]{geometry}
\usepackage[T1]{fontenc}
\usepackage{lmodern,microtype}

\usepackage[dvipsnames]{xcolor}
\usepackage[
 colorlinks=true,
 linkcolor=MidnightBlue,
 citecolor=MidnightBlue,
 urlcolor=MidnightBlue
]{hyperref}

\numberwithin{equation}{section}
\newcommand{\xx}{\mathbf{x}}
\newcommand{\yy}{\mathbf{y}}
\renewcommand{\aa}{\mathbf{a}}
\newcommand{\zz}{\mathbf{z}}
\DeclareMathOperator{\rank}{rank}
\DeclareMathOperator{\Br}{Br}
\renewcommand{\leq}{\leqslant}
\renewcommand{\geq}{\geqslant}
\newcommand{\ZZ}{\mathbb{Z}}
\newcommand{\QQ}{\mathbb{Q}}
\newcommand{\CC}{\mathbb{C}}

\newcommand{\FF}{\mathbb{F}}
\newcommand{\A}{\mathbb{A}}
\newcommand{\PP}{\mathbb{P}}
\newcommand{\Sha}{\mathop{\mathrm{\amalg}}}
\newtheorem{theorem}{Theorem}[section]
\newtheorem{corollary}[theorem]{Corollary}
\newtheorem{lemma}[theorem]{Lemma}
\newtheorem{proposition}[theorem]{Proposition}

\theoremstyle{remark}
\newtheorem{remark}[theorem]{Remark}

\title[Rational lines on pairs  of quadrics and cubic hypersurfaces]{Rational lines on pairs  of quadrics\\ and cubic hypersurfaces}
\author{Tim Browning}

\address{ISTA\\
Am Campus 1\\
3400 Klosterneuburg\\
Austria}
\email{tdb@ist.ac.at}

\address{Department of Mathematics\\ Imperial College London\\
SW7 2AZ\\ 
United Kingdom}
\email{tim.browning@imperial.ac.uk}

		\author{Stephanie Chan}
\address{Department of Mathematics,
University College London,
Gower Street,
London WC1E 6BT, United Kingdom}
\email{stephanie.chan@ucl.ac.uk}
		\author{Nick Rome}
        \address{Max Planck Institute for Mathematics, Vivatsgasse 7, 53111 Bonn, Germany
        }
        \email{nrome@mpim-bonn.mpg.de}

 \subjclass[2020]{11D72 (11E12, 14G05, 14G20)}
\date{}

\begin{document}
\begin{abstract}
We study rational lines on smooth intersections of two quadrics and on
smooth cubic hypersurfaces over number fields. For an intersection of
two quadrics in $n$ variables, we prove the Hasse principle and weak
approximation for the variety of lines when $n\geq11$, and more generally for the variety of $r$-planes when $n\geq3r+8$. In this range an $r$-plane exists over the ground field if and only if one exists over every
real completion. For the variety of lines, we also prove the Hasse principle and weak approximation when $n\geq 8$, conditional on the finiteness
of certain Tate--Shafarevich groups. For a smooth cubic hypersurface in $n$ variables, we prove that a line
exists over the ground field when $n\geq23$. 
\end{abstract}
\maketitle
\thispagestyle{empty}

 \setcounter{tocdepth}{1}
 \tableofcontents

\section{Introduction}

Let $k$ be a number field and let $X\subset\PP^{n-1}_k$ be a smooth
projective variety. For an integer $r\geq1$, we write $F_r(X)$ for the
Fano variety parametrising $r$-dimensional projective linear subspaces
$\PP^r\subset X$. Thus a point of $F_r(X)(k)$ is an $(r+1)$-dimensional
vector subspace of $k^n$ on which the defining equations of $X$ vanish
identically. When $r=1$, $F_1(X)$ is the variety of lines on $X$.
This paper treats two natural cases, namely smooth complete
intersections of two quadrics and smooth cubic hypersurfaces. In each
case we use a rational point of $X$ to reduce the problem to a system of
two equations in fewer variables.
The following is our first main result. 

\begin{theorem}\label{thm:quadrics}
Let $X\subset\PP^{n-1}_k$ be a smooth complete intersection of two
quadrics, with $n\geq11$. Then
\[
 F_1(X)(k)\ne\varnothing
 \quad\Longleftrightarrow\quad
 F_1(X)(k_v)\ne\varnothing\quad\text{for every real place $v$ of $k$}.
\]
If these conditions hold, then $F_1(X)(k)$ is dense in $F_1(X)(\A_k)$.
In particular, $F_1(X)$ satisfies the Hasse principle and weak
approximation.
\end{theorem}

The real place condition is clearly vacuous when $k$ is totally imaginary.
Hence every smooth intersection of two quadrics in at least $11$
variables over such a field contains a line. 
In general, 
a concrete way to check the 
real place condition is provided by Ji and Suzuki
\cite[Lemma~5.6(1)]{ji}. 
Over $\QQ$ and for $n\geq11$, it follows from Theorem~\ref{thm:quadrics} that the existence of a rational line
is equivalent to the existence of a real line. 
Consider the smooth intersection of two quadrics defined by 
\[
 Q_0=x_1^2-\sum_{j=2}^n x_j^2,\qquad
 Q_1=16x_2^2-9x_3^2-\sum_{j=4}^n(j+10)x_j^2.
\]
This contains the rational point 
$[5:3:4:0:\cdots:0]$, but $Q_0$ has signature $(1,n-1)$ and cannot contain a
two-dimensional totally isotropic subspace.  This  therefore 
gives an example of  a smooth intersection over $\QQ$ with a rational point and no real line. 

In fact, we can iterate our
fibration argument by induction on $r$ to handle $r$-planes on $X$ for
every $r\geq1$.

\begin{corollary}\label{cor:rplanes}
Let $r\geq1$ and let $X\subset\PP^{n-1}_k$ be a smooth complete
intersection of two quadrics, with $n\geq3r+8$. Then
\[
 F_r(X)(k)\ne\varnothing
 \quad\Longleftrightarrow\quad
 F_r(X)(k_v)\ne\varnothing\quad\text{for every real place $v$ of $k$}.
\]
If these conditions hold, then $F_r(X)(k)$ is dense in $F_r(X)(\A_k)$.
In particular, $F_r(X)$ satisfies the Hasse principle and weak
approximation.
\end{corollary}

Over a totally imaginary number field, Dietmann and Wooley
\cite[Equation (2.11)]{dw} showed that any pair of quadratic forms in
$n\geq 3r+9$ variables vanishes on an  $r$-plane defined over the number field. 
For smooth intersections, Corollary~\ref{cor:rplanes} saves one variable
and also handles arbitrary number fields.

The arithmetic of points on intersections of two quadrics is well
understood in sufficiently many variables. The work of
Colliot-Th\'el\`ene, Sansuc and Swinnerton-Dyer
\cite{css1,css2}, followed by Heath-Brown \cite{hb}, establishes the
Hasse principle for smooth intersections in $\PP^{n-1}$ when $n\geq 8$.
These results form the arithmetic input to Theorem~\ref{thm:quadrics}.
In \cite[Theorem~3.8]{moly}, 
Molyakov
has extended this to smooth proper models of geometrically integral,
non-conical intersections in $\PP^7$, but we shall not need this in our work. 
Our work is inspired by work of Ji and
Suzuki \cite{ji}. Their Corollaries~3.8 and~3.10 show, for $n\geq2r+6$,
that $\Br(k)\to\Br(F_r(X))$ is surjective, and that weak approximation
holds as soon as $F_r(X)(k)$ is nonempty. 
Colliot-Th\'el\`ene's conjecture
\cite{bud}
on the Brauer--Manin obstruction therefore predicts the Hasse principle
already for $n\geq2r+6$. 
We are able to confirm this in the case $r=1$ conditionally, assuming the finiteness of certain Tate--Shafarevich groups.

\begin{theorem}\label{thm:shaquadrics}
    Let $n \geq 8$ and let $X \subset \PP_k^{n-1}$ be a smooth complete intersection of two quadrics. 
    Assume that
\[
\Sha(\operatorname{Jac}(C)/L)[2^\infty]
\]
is finite for every finite extension $L/k$ and every genus
$2$ curve $C/L$.
    Then $F_1(X)$ satisfies the Hasse principle and weak approximation.
\end{theorem}

The condition in the above theorem is used to control the arithmetic of the intersections of two quadrics arising as fibres in our argument.
In particular, work of Morgan--Skorobogatov~\cite{ms} shows that smooth complete intersections of two quadrics in $\PP^{n-1}$ for $n \geq 6$ satisfy the Hasse principle under this assumption.

In smaller dimension the
Hasse principle can fail for $F_1(X)$. 
When $n=6$, the variety of lines is a torsor under the Jacobian of
a genus $2$ curve \cite[Theorem~2.2(2)]{ji}. 
Moreover, 
Frei and Ji
\cite{freiji} give an explicit smooth intersection
in $\PP^5$ over $\QQ$ which contains a line over every completion
of $\QQ$ but contains no rational line.

\medskip

Our second main result concerns cubic hypersurfaces.

\begin{theorem}\label{thm:cubic}
Let $X\subset\PP^{n-1}_k$ be a smooth cubic hypersurface, with
$n\geq23$. Then $X$ contains a $k$-line.
\end{theorem}

Equivalently, if $C\in k[x_1,\ldots,x_{n}]$ is a nonsingular cubic
form in at least $23$ variables, there are linearly independent vectors
$\xx,\yy\in k^n$ such that
$C(s\xx+t\yy)=0$ for all $s,t\in k$.
Note that there are no local hypotheses in this result.

Conditionally, the Hasse principle holds on $F_1(X)$ once $n \geq 7$. Indeed, $F_1(X)$ is a smooth projective Fano variety for $n \geq 7$ by 
\cite[Chapter~2, Lemma~3.1 and Proposition~3.4]{huy}.
Thus, Colliot-Th{\'e}l{\`e}ne's conjecture \cite{bud} predicts that the Brauer--Manin obstruction is the only obstruction to the Hasse principle. 
Over $\CC$, Kodaira vanishing gives $H^2(F_1(X)_{\CC},\mathcal O_{F_1(X)})=0$,
and $H^3(F_1(X)_{\CC},\ZZ)$ is torsion-free by
\cite[Chapter~2, Exercise~4.13]{huy}. The exponential sequence
therefore gives $\Br(F_1(X)_{\CC})=0$, whence $\Br(\overline{F_1(X)})=0$
by comparison. Moreover, $\operatorname{Pic}(\overline{F_1(X)})\cong\ZZ$
by \cite[Corollary~3.5 and Remark~3.6(5)]{DeMa}.
The Galois action is trivial since it fixes the Pl\"ucker class.
Thus $H^1(k,\operatorname{Pic}(\overline{F_1(X)}))=0$ and the
Hochschild--Serre sequence implies that $\Br(k)\to\Br(F_1(X))$ is
surjective; it follows that  $F_1(X)(\A_k)^{\Br}=F_1(X)(\A_k)$.

An aspect that is less understood, even conjecturally, is the least number of variables required to guarantee the existence of a line over the non-archimedean completions of $k$. 
Over $\QQ_p$, Brandes and Dietmann \cite[Theorems~1.3 and~1.4]{BD}
show that every cubic form in at least $22$ variables vanishes on a
line, and that $16$ variables suffice for  $p\geq19$.
For each $p\in\{2,3,5\}$ they also give a smooth cubic in $11$
variables with no $\QQ_p$-line.

Wooley \cite[Theorem~2(b)]{trevor} proved that every cubic
form in at least $37$ variables over a number field
vanishes on a rational line.
Over $\QQ$, Brandes and Dietmann
\cite[Theorem~1.1]{BD} proved that $31$ variables suffice
for a nonsingular cubic form.
Brandes, Dietmann
and Leep \cite[Theorem~2]{BDL} subsequently obtained bounds of $35$
variables over any number field, $33$ over an imaginary quadratic
field, and $31$ over $\QQ$, without assuming nonsingularity. 
Their
argument constructs a line over a suitable quadratic extension and
then descends it. 
As Brandes, Dietmann and Leep observe after \cite[Theorem~2]{BDL},
$23$ is the limiting number of variables for their method using a
general quadratic form. We attain this bound for smooth cubics by
choosing the linear space on the quadric so that the cubic restricted
to it is nonsingular. We can then apply the ten variable theorem of
Browning and Vishe \cite[Theorem~1.1]{BV}.

\medskip

We proceed by describing the common starting point. Let
\begin{equation}\label{eq:incidence}
 I=\{(a,\ell)\in X\times F_1(X):a\in\ell\},
 \qquad \pi:I\longrightarrow X.
\end{equation}
The fibre over $a$ parametrises the lines on $X$ through $a$.
If $X=\{Q_0=Q_1=0\}\subset \PP^{n-1}$ 
is a smooth intersection of two quadrics  over a number field $k$ 
and $a=[\aa]$, write $B_i$ for the symmetric
bilinear form with $Q_i(\xx)=B_i(\xx,\xx)$. The conditions on a second
vector are 
\[
 B_0(\aa,\yy)=B_1(\aa,\yy)=Q_0(\yy)=Q_1(\yy)=0,
 \qquad \yy\notin k\aa.
\]
The linear equations leave $n-2$ variables. Removing the direction
$k\aa$ leaves an intersection of two quadrics $Y_a\subset\PP^{n-4}$, which will be smooth for general $a$. Thus $n=11$ is precisely
the point at which Heath-Brown's theorem applies to the fibre. To
choose $a$ with an everywhere locally soluble fibre, we prove that
outside a fixed finite set of places, \emph{every} fibre has a smooth
local point. Weak approximation on $X$ then handles the remaining
places. 
This will be carried out in Section~\ref{sec:pairs}.

For a cubic hypersurface, choose coordinates at $a$ so that
\[
 a=[1:0:\cdots:0],\qquad T_aX=\{x_2=0\}.
\]
Writing $\zz=(x_3,\ldots,x_{n})$, the tangent section has equation
\[
 C(x_1,0,\zz)=x_1q(\zz)+c(\zz),
\]
and the lines through $a$ are parametrised by $q=c=0$ in
$\PP^{n-3}$. We choose $a$ so that $q$ is nonsingular and this
intersection is smooth. Theorem~3 of \cite{BDL} provides a quadratic
extension over which a quadratic form in $2d+1$ variables contains a
$d$-dimensional vector space of zeros. We show that this space can be
chosen such that the restriction of the cubic to this space is nonsingular. Taking $d=10$ gives
$n-2\geq21$, whence  $n\geq23$.
The argument for cubics is explained in 
Section~\ref{sec:cubic}. 


\subsection*{Acknowledgements}
NR would like to thank Professor Colliot-Th{\'e}l{\`e}ne for pointing him to the reference \cite{DeMa} in regards to the Brauer group of $F_1(X)$. During part of this work, NR was funded by the FWF project ESP 441-NBL.

\subsection*{AI statement}

The structure of the proofs and most of the ideas in this paper are due to the authors. 
ChatGPT 6.0 suggested 
Proposition \ref{prop:uniform} and its proof, 
as a means of establishing 
a sufficiently uniform version of local solubility for the intersection of two quadrics $Y_a\subset \PP^{n-4}$ that emerges in the proof of Theorem  \ref{thm:quadrics}. 
Moreover, ChatGPT 6.0   suggested the crucial nonsingularity 
result in Lemma \ref{lem:restriction} that is needed to 
apply the ten variable theorem for cubics over number fields. 
Gemini 4 suggested generalising from lines to $r$-planes on pairs of quadrics, and helped to provide shorter proofs for Lemmas~\ref{lem:rank}, \ref{lem:finitefield} and \ref{lem:base-point}.

\section{Pairs of quadrics}\label{sec:pairs}

\subsection{Lines through a point on two quadrics}

In this section let
\begin{equation}\label{eq:X}
 X=\{Q_0=Q_1=0\}\subset\PP^{n-1}_K
\end{equation}
be a smooth complete intersection over a field $K$ of characteristic
different from $2$. Write $B_i$ for the associated symmetric bilinear
form, with $Q_i(\xx)=B_i(\xx,\xx)$. All
geometric assertions in this section may be checked over an algebraic
closure $\overline K$.

\begin{lemma}\label{lem:pencil}
If $X$ in \eqref{eq:X} is smooth, then
$\det(sQ_0+tQ_1)$ is a nonzero binary form of degree $n$ with distinct
geometric roots. Every member of the geometric pencil has rank at least
$n-1$.
\end{lemma}

\begin{proof}
This is \cite[Lemma~3.1]{hb}.
\end{proof}

For $a=[\aa]\in X(K)$, put
\begin{equation}\label{eq:Wa}
 W_a=\{\yy\in K^n:B_0(\aa,\yy)=B_1(\aa,\yy)=0\},
 \qquad E_a=W_a/K\aa.
\end{equation}
Since $X$ is a smooth complete intersection, the two gradients $\nabla Q_0(\mathbf a)$ and 
$\nabla Q_1(\mathbf a)$ are linearly independent. Thus the equations $2B_i(\mathbf a,\mathbf y)=\nabla Q_i(\mathbf a)\cdot\mathbf y=0$ cut out a space $W_a$ of dimension $n-2$.
Moreover, $B_i(\mathbf a,\mathbf a)=Q_i(\mathbf a)=0$ for $i=0,1$, so $W_a$ contains the line 
$K\mathbf a$. Taking the quotient by this line gives $\dim E_a=n-3$.
For $\yy\in W_a$ and $c\in K$, we clearly have 
\[
 Q_i(\yy+c\aa)=Q_i(\yy)+2cB_i(\aa,\yy)+c^2Q_i(\aa)=Q_i(\yy),
\]
for $i=0,1$. It follows that  the restrictions induce quadratic forms $R_{0,a},R_{1,a}$
on $E_a$. Define
\begin{equation}\label{eq:Ya}
 Y_a=\{R_{0,a}=R_{1,a}=0\}\subset\PP(E_a).
\end{equation}
Choosing a different nonzero representative of $a$ changes neither
$W_a$ nor this construction.

\begin{lemma}\label{lem:rank}
For every geometric point $a\in X$ and every $[s:t]\in\PP^1$ we have
\[
 \rank(sR_{0,a}+tR_{1,a})\geq n-5,
\]
and $\det(sR_{0,a}+tR_{1,a})$ is not identically zero.
If $n\geq8$, then $Y_a$ is a geometrically integral complete
intersection of dimension $n-6$.
If $n=8$ and $\operatorname{char}K=0$, then at the generic point
$\eta$ of $X$ the polynomial
$\Delta_\eta(T)=\det(R_{1,\eta}-TR_{0,\eta})$ is irreducible
of degree $5$ over $\overline K(X)$.
\end{lemma}

\begin{proof}
We work over $\overline K$. By Lemma~\ref{lem:pencil}, every member
of the original pencil has rank at least $n-1$. Restricting to
$W_a$, which has codimension two, lowers the rank by at most four
\cite[Lemma~1.8]{css1}. Passing to $E_a$ leaves the rank unchanged,
since $\aa$ is already in the kernel. This gives the bound $n-5$.

By Lemma~\ref{lem:pencil}, the pencil has only finitely many
singular members over $\overline K$. We may therefore choose
two distinct nonsingular members as generators, and denote them
again by $Q_0,Q_1$. By \cite[Proposition~2.1]{hbp}, we may
simultaneously diagonalise them and rescale the coordinates. 
In these new coordinates, we have 
\[
 Q_0(\yy)=\sum_{i=1}^n y_i^2,\qquad
 Q_1(\yy)=\sum_{i=1}^n d_i y_i^2,
\]
where the $d_i$ are pairwise distinct.
To prove that the determinant is not identically zero, it suffices
to find one member of the pencil that is nonsingular on $E_a$.
Since $E_a=W_a/\langle\aa\rangle$, this amounts to finding a member
whose kernel on $W_a$ is exactly $\langle\aa\rangle$.
Since $\aa\ne \mathbf 0$, we may relabel the coordinates so that $a_n\ne0$.
We claim that $Q=Q_1-d_nQ_0$ has this property.
We clearly have 
\[
 Q(\yy)
       =\sum_{i=1}^{n-1}(d_i-d_n)y_i^2.
\]
Since the $d_i$ are distinct, this is a nonsingular quadratic
form in $y_1,\ldots,y_{n-1}$.
Since $a_n\ne0$, the equation $B_0(\aa,\yy)=0$ determines $y_n$.
Subtracting $d_nB_0(\aa,\yy)$ from $B_1(\aa,\yy)$ shows that
the first $n-1$ coordinates of $\yy\in W_a$ satisfy 
\[
 \sum_{i<n}(d_i-d_n)a_i y_i=0.
\]
Suppose that $\yy$ lies in the kernel of the restriction of $Q$
to $W_a$. Then the linear form
$\sum_{i<n}(d_i-d_n)y_i z_i$ vanishes whenever
$\sum_{i<n}(d_i-d_n)a_i z_i=0$, so the former is a scalar
multiple of the latter. Since $d_i\ne d_n$ for $i<n$, comparison
of coefficients gives $y_i=\lambda a_i$ for $i<n$. Finally,
\[
 a_n y_n=-\sum_{i<n}a_i y_i
        =-\lambda\sum_{i<n}a_i^2=\lambda a_n^2,
\]
since $Q_0(\aa)=0$. Thus $\yy=\lambda\aa$. As $\aa$ already
lies in the kernel, the kernel is exactly $\langle\aa\rangle$,
and $Q$ induces a nonsingular form on $E_a$.

If $n\geq8$, every member of the induced pencil has rank at least
$n-5\geq3$, and one has rank $n-3\geq5$. 
The two generators have no common factor. Indeed, a common linear
factor would force each to have rank at most two, while a common
quadratic factor would make them proportional, contrary to the
rank bound for their linear combinations. We may therefore apply
\cite[Lemma~1.11]{css1}
to deduce that $Y_a$ is geometrically
integral of codimension two in $\PP(E_a)$. It is therefore a
complete intersection of dimension $n-6$.

Finally, suppose that $n=8$ and $\operatorname{char}K=0$.
By Jeong \cite[Proposition~3.9(i)]{jeong},
the map taking $a$ to
its binary determinant sends $[a_1:\cdots:a_8]$ to
$[a_1^2:\cdots:a_8^2]$, followed by a linear isomorphism from
the image onto the projective space of binary quintics.
(Jeong works over $\CC$, but  the same formula holds in
characteristic zero by embedding a finitely generated field of
definition into $\CC$.)
Since this map is dominant, the generic point $\eta$ of $X$
maps to the generic point of the space of binary quintics.
Dehomogenising and dividing by the leading coefficient, we
therefore obtain
\[
 T^5+c_4T^4+c_3T^3+c_2T^2+c_1T+c_0,
\]
where $c_0,\ldots,c_4$ are algebraically independent over $\overline K$.
This is the generic monic quintic, which is irreducible over
$F=\overline K(c_0,\ldots,c_4)$.
 Since the squared coordinates $[a_1^2:\cdots:a_8^2]$ can be
recovered from the coefficients of the quintic by a linear
isomorphism, each ratio $a_i^2/a_8^2$ belongs to $F$.
To obtain $\overline K(X)$ we need only adjoin the ratios
$a_i/a_8$ for $1\leq i\leq7$. Each is a square root of an
element of $F$, so $[\overline K(X):F]$ is a power of two.
Put $E=\overline K(X)$ and let $\alpha$ be a root of the quintic.
Since $[F(\alpha):F]=5$, the tower law gives
\[
 5\mid [E(\alpha):F]=[E(\alpha):E][E:F].
\]
As $[E:F]$ is a power of two, $5$ divides $[E(\alpha):E]$.
This degree is at most five, since $\alpha$ satisfies the quintic,
so it is exactly five. Thus the quintic is irreducible over $E$,
and so is its homogenisation, the binary determinant.
For the original generators, put
\[
 D(U,V)=\det(VR_{1,\eta}-UR_{0,\eta}),
\]
so that $D(T,1)=\Delta_\eta(T)$. Changing generators amounts to
an invertible linear change in $U,V$, so $D$ is irreducible over
$\overline K(X)$. In particular, $V$ does not divide $D$, so its
$U^5$ coefficient is nonzero. Thus $\Delta_\eta$ has degree five,
and any factorisation of it would homogenise to a factorisation
of $D$. This proves the lemma.
\end{proof}

Recall the incidence variety $I$ and projection $\pi:I\to X$ from \eqref{eq:incidence}.
The other projection $I\to F_1(X)$ is a projective line bundle.
The fibre $\pi^{-1}(a)$ is naturally isomorphic to $Y_a$. Indeed,  the point
$[\yy]\in Y_a$ corresponds to the line
$\PP(K\aa+K\yy)$. 
This is an isomorphism of schemes, as is seen by
identifying lines through $a$ with one-dimensional subspaces of
$K^n/K\aa$ and then imposing the two linear and two quadratic equations
for a line to lie on $X$.
More generally, for any integer $r\geq1$, consider the incidence variety
\begin{equation}\label{eq:incidence-r}
 I_r=\{(a,\Lambda)\in X\times F_r(X):a\in\Lambda\},
 \qquad \pi_r:I_r\longrightarrow X
\end{equation}
so that $I_1=I$ and $\pi_1=\pi$.
Sending $\Lambda=\PP(V)$ to $\PP(V/K\aa)\subset\PP(E_a)$ identifies the fibre $\pi_r^{-1}(a)$ with the variety $F_{r-1}(Y_a)$ of $(r-1)$-planes on $Y_a$, with the convention $F_0(Y_a)=Y_a$.

\begin{lemma}\label{lem:generic}
Suppose that $\operatorname{char}K=0$ and $n\geq \max\{8,2r+5\}$ for an integer $r\geq 1$. The varieties
$F_r(X)$ and $I_r$ are smooth and geometrically integral. There is a
nonempty Zariski open subset $U\subset X$, defined over $K$, over which $\pi_r$ is smooth and $Y_a$ is a smooth geometrically integral intersection of two
quadrics for every $a\in U$.
\end{lemma}

\begin{proof}
For a smooth complete intersection of two quadrics $V\subset\PP^N$, the
scheme $F_s(V)$ is smooth, projective and nonempty if $N\geq 2s+2$,
and geometrically connected if $N\geq 2s+3$ by
\cite[Lemma~2.1]{ji}. Thus $F_r(X)$ is smooth and geometrically
integral, and the same holds for the $\PP^r$-bundle $I_r$.

By \cite[Proposition~2.14(f)]{ctretour}, there is a nonempty open
subset $U\subset X$ over which $Y_a$ is a smooth geometrically
integral complete intersection of two quadrics. For each geometric
point $a\in U$, the fibre of $\pi_r$ is $F_{r-1}(Y_a)$.
This is nonempty by the first assertion above, since
$n-4\geq 2(r-1)+2$. Hence $\pi_r$ is dominant.
Shrinking $U$ once more, generic smoothness in characteristic zero
gives the result.
\end{proof}

\begin{remark}
The intersection $X\cap\PP(W_a)$ is the projective cone over $Y_a$
with vertex $a$. Returning to  \eqref{eq:Wa} and \eqref{eq:Ya}, 
we have 
$\dim W_a=n-2$ and $\dim E_a=n-3$, with 
 $Y_a\subset\PP^{n-4}$.
\end{remark}

\subsection{Uniform local solubility}

The following finite field estimate is the crucial uniform input needed in the
fibration argument to come.

\begin{lemma}\label{lem:finitefield}
Suppose that $n\geq 8$ and $K=\FF_q$, where $q$ is odd and
sufficiently large in terms of $n$. Then $Y_a$ has a smooth
$\FF_q$-point for every $a\in X(\FF_q)$.
\end{lemma}

\begin{proof}
Consider the family over $\ZZ$ of common zero schemes of two quadrics
in $\PP^{n-4}$, with the coefficients of both quadrics as parameters.
Choosing coordinates on $E_a$ identifies $Y_a$ with a fibre of this
family. By Lemma~\ref{lem:rank}, this fibre is geometrically integral.
Thus
the uniform Lang--Weil theorem 
 \cite[Theorem~7.7.1(iv)]{poonen} gives a smooth $\FF_q$-point
for $q$ sufficiently large in terms of the family, hence in terms
of $n$ alone.
\end{proof}

\begin{proposition}\label{prop:uniform}
Let $k$ be a number field and suppose that $n\geq8$. There is a finite
set $S_0$ of places of $k$, containing every archimedean place, depending only on the equations of $X$, such
that for every place $v\notin S_0$ and every $a\in X(k_v)$,
the variety $Y_a$ has a smooth $k_v$-point.
\end{proposition}

\begin{proof}
Choose $S_0$ to contain the archimedean places, the places above $2$, the places where the given equations are not integral or do not define a smooth complete intersection after reduction, and the places whose residue field $\FF_q$ fails the bound of Lemma~\ref{lem:finitefield}. For $v\notin S_0$ with valuation ring $\mathcal O_v$, a primitive vector $\aa\in\mathcal O_v^n$ representing $a\in X(k_v)$ reduces to a point $\overline a$ on the smooth reduction of $X$, so $\mathcal E_a=\{\yy\in\mathcal O_v^n:B_0(\aa,\yy)=B_1(\aa,\yy)=0\}/\mathcal O_v\aa$ is free of rank $n-3$ and reduces to $E_{\overline a}$. By Lemma~\ref{lem:finitefield}, $Y_{\overline a}$ has a smooth $\FF_q$-point, which lifts to a smooth $k_v$-point of $Y_a$ by Hensel's lemma \cite[Theorem~3.5.63]{poonen}.
\end{proof}

\begin{remark}
Note that the set $S_0$ is chosen \emph{before} choosing $a$; proving that each fixed smooth $Y_a$ has local points at
almost all places would only give a finite exceptional set depending on
$a$, and would not justify the later use of weak approximation.
\end{remark}

\subsection{Final steps}

We now record the results on the  arithmetic of intersections of two quadrics that we will need in our work. 

\begin{proposition}\label{prop:points}
Let $V\subset\PP^N$ be a smooth complete intersection of two quadrics
over a number field $k$.
\begin{enumerate}
\item If $N\geq7$, then $V$ satisfies the Hasse principle.
\item If $r\geq 0$, $N\geq 2r+5$, and $F_r(V)(k)\ne\varnothing$, then $F_r(V)$ satisfies weak approximation. 
\end{enumerate}
\end{proposition}

\begin{proof}
The first assertion for $N=7$ is \cite[Theorem~1]{hb}; the case
$N\geq8$ is due to Colliot-Th\'el\`ene--Sansuc--Swinnerton-Dyer
\cite[Theorem~A(i)(a)]{css1}. The second is \cite[Corollary~3.10]{ji}. For $r=0$, it follows from \cite[Theorem~A(ii)]{css1}.
\end{proof}

We will use the following standard consequences of the implicit
function theorem over local fields; see
\cite[Proposition~3.5.73(ii) and Remark~3.5.76]{poonen}.

\begin{lemma}\label{lem:local}
Let $K$ be a completion of a number field. If $V$ is a smooth
geometrically integral $K$-variety and $V^\circ\subset V$ is a nonempty
Zariski open subset, then $V^\circ(K)$ is dense in $V(K)$.
If $f:V\to W$ is a morphism of smooth $K$-varieties that is smooth at
$x\in V(K)$, then the image under $f$ of every sufficiently small
neighbourhood of $x$ contains a neighbourhood of $f(x)$.
\end{lemma}

Let us return now to the setting of our smooth intersection of two quadrics $X\subset \PP_k^{n-1}$
over a number field $k$. 
Recall the definition \eqref{eq:incidence-r} of the incidence variety $I_r$.
Armed with the above  results, we  now describe how to choose our  base point $a$. 

\begin{proposition}\label{prop:selection}
Suppose that $r\geq1$, $n\geq\max\{8,2r+5\}$ and
$F_r(X)(k_v)\ne\varnothing$ for every place $v$ of $k$.
Let $U$ be as in Lemma~\ref{lem:generic}.
There exists $a\in U(k)$ such that $F_{r-1}(Y_a)(k_v)\ne\varnothing$
at every archimedean place $v$. When $r=1$, we may take $Y_a$ to have
a point over every completion of $k$. If also $n=8$, we may take
$\Delta_a(T)$ to be irreducible over $k$.
\end{proposition}

\begin{proof}
Each local $r$-plane contains a local point. Thus $X$ is everywhere
locally soluble, and Proposition~\ref{prop:points} gives a $k$-point
and weak approximation on $X$. 
When $(r,n)=(1,8)$, Lemma~\ref{lem:rank} allows us to shrink $U$ so
that $R_{0,a}$ is nonsingular throughout $U$. We use the monic
representative
\[
 \frac{\det(TR_{0,a}-R_{1,a})}{\det R_{0,a}}
\]
of $\Delta_a(T)$. Its coefficients are regular on $U$, since a
change of basis multiplies both determinants by the same factor.

Let $S_0$ be the finite set in Proposition~\ref{prop:uniform}, and put
$I_{r,U}=\pi_r^{-1}(U)$. For $v\in S_0$, a local $r$-plane together
with a point on it gives a point of $I_r(k_v)$. Since $I_r$ is smooth
and geometrically integral, Lemma~\ref{lem:local} allows us to move
this point into $I_{r,U}(k_v)$. The map $\pi_r$ is smooth there, so
\[
 A_v=\pi_r(I_{r,U}(k_v))
\]
is a nonempty open subset of $X(k_v)$.

When $(r,n)=(1,8)$ the generic polynomial is irreducible by
Lemma~\ref{lem:rank}. Hence the specialisations in $U(k)$ for which it
is reducible lie in a thin subset $Z\subset X(k)$ by
\cite[Proposition~3.3.5]{serre}. In the other cases put $Z=\varnothing$.
Since $X$ is projective, the closure of $Z$ in $\prod_v X(k_v)$ has
empty interior by \cite[Corollary~3.5.4]{serre}. By weak approximation
we can therefore choose $a\in X(k)\setminus Z$ lying in $A_v$ for
every $v\in S_0$. In particular, $a\in U(k)$ and the fibre
$F_{r-1}(Y_a)$ has a $k_v$-point for every $v\in S_0$. When $r=1$,
Proposition~\ref{prop:uniform} supplies a point at every remaining place.
\end{proof}

\begin{proposition}\label{prop:transfer}
For $n \geq 11$, $F_1(X)$ satisfies the Hasse principle and weak approximation. Moreover,
if $n \geq 8$ and the hypothesis of Theorem~\ref{thm:shaquadrics} holds, then $F_1(X)$
satisfies the Hasse principle and weak approximation.
\end{proposition}

\begin{proof}
Suppose that $F_1(X)$ has a point over every completion of $k$.
By Proposition~\ref{prop:selection} with $r=1$, there is a point
$a\in X(k)$ for which $Y_a\subset\PP^{n-4}_k$ is a smooth,
everywhere locally soluble intersection of two quadrics.
If $n\geq11$, then $n-4\geq7$, so $Y_a(k)\ne\varnothing$ by
Proposition~\ref{prop:points}(1). For $8\leq n\leq10$, assume the
finiteness hypothesis of Theorem~\ref{thm:shaquadrics}.
If $n=9$ or $10$, then $Y_a$ lies in $\PP^5$ or $\PP^6$ and
\cite[Theorem~B]{ms} gives a $k$-point. If $n=8$, then $Y_a$ is a
quartic del Pezzo surface in $\PP^4$. In this case
Proposition~\ref{prop:selection} allows us to choose $a$ so that
$\Delta_a(T)$ is irreducible, and \cite[Theorem~A]{ms} gives a
$k$-point.
In each case a $k$-point of $Y_a$ gives a $k$-line on $X$ through $a$.
Weak approximation on $F_1(X)$ now follows from
Proposition~\ref{prop:points}(2).
\end{proof}

Proposition~\ref{prop:transfer} proves 
Theorem~\ref{thm:shaquadrics}.
For $n\geq11$, it gives the Hasse principle and weak
approximation for $F_1(X)$ unconditionally. The real place criterion follows from
the next local result, which builds on work of 
 Parimala--Suresh \cite[Theorem~4.6]{PS} and 
 Leep \cite[Theorem~3.4]{Leep} on the $u$-invariant of a $p$-adic function field.

\begin{lemma}\label{lem:padic}
Let $p$ be a prime, let $K$ be a finite extension of $\QQ_p$, and let $r\geq 1$.
If $n\geq 3r+8$, then every smooth intersection of two quadrics in $\PP^{n-1}_K$
contains an $r$-dimensional projective linear subspace over $K$.
\end{lemma}

\begin{proof}
This follows from \cite[Remarks~3.9(3)]{ji}, taking
$N=n-1$; the argument in the proof of
\cite[Corollary~3.5(2)]{ji} applies over every finite extension
of $\QQ_p$.
\end{proof}

\begin{proof}[Proof of Theorem~\ref{thm:quadrics} and Corollary~\ref{cor:rplanes}]
Since Theorem~\ref{thm:quadrics} is the case $r=1$ of Corollary~\ref{cor:rplanes} and necessity is clear, we prove sufficiency by induction on $r\geq 1$. Suppose that $n\geq 3r+8$ and $F_r(X)(k_v)\ne\varnothing$ for every real place $v$. Lemma~\ref{lem:padic} supplies $r$-planes at every finite place, and $r$-planes also exist over $\CC$ by the geometric nonemptiness of $F_r(X)$ in Lemma~\ref{lem:generic}. Thus $F_r(X)$ is everywhere locally soluble. When $r=1$, Proposition~\ref{prop:transfer} gives a rational line and weak approximation. When $r\geq 2$, Proposition~\ref{prop:selection} gives $a\in U(k)$ such that the smooth complete intersection $Y_a\subset\PP^{n-4}_k$ in $n-3\geq 3(r-1)+8$ variables satisfies $F_{r-1}(Y_a)(k_v)\ne\varnothing$ for every real place $v$. By the induction hypothesis, $F_{r-1}(Y_a)(k)\ne\varnothing$, which gives $F_r(X)(k)\ne\varnothing$ and weak approximation by Proposition~\ref{prop:points}(2).
\end{proof}

\section{Cubic hypersurfaces}\label{sec:cubic}

\subsection{Lines through a point on a cubic hypersurface}

Let $K$ be a field of characteristic zero and let
$X=\{C=0\}\subset\PP^{n-1}_K$ be a smooth cubic hypersurface, with
$n\geq5$.
At a point $a\in X(K)$, choose coordinates in which
\[
 a=[1:0:\cdots:0],
 \qquad T_{a}X=\{x_2=0\}.
\]
Writing $\zz=(x_3,\ldots,x_{n})$, we then have
\begin{equation}\label{eq:tangent-section}
 C(x_1,0,\zz)=x_1q(\zz)+c(\zz),
\end{equation}
where $q$ and $c$ are quadratic and cubic forms, respectively.
A line through $a$
is contained in $X$ precisely when its direction gives a point on
\[
 Y_{a}=\{q=c=0\}\subset\PP^{n-3}.
\]

\begin{lemma}\label{lem:base-point}
There is a nonempty open subset $X^\circ\subset X$ defined over $K$ such that, for
$a\in X^\circ$, the quadratic form $q$ in
\eqref{eq:tangent-section} is nonsingular and $Y_a$ is a smooth complete intersection.
\end{lemma}

\begin{proof}
Work over $\overline K$. By \cite[Chapter~1, Remark~5.20]{huy},
for general $a$ the tangent hyperplane section $X\cap T_aX$ has just one
singular point, an ordinary double point at $a$. Thus the quadratic
part $q$ in \eqref{eq:tangent-section} is nonsingular.
Suppose that $[\zz]\in Y_a$ is a point where the gradients of $q$
and $c$ are dependent. Since $q$ is nonsingular, we have
$\nabla c(\zz)=\lambda\nabla q(\zz)$ for some $\lambda$.
Then $[-\lambda:0:\zz]$ is another singular point of the tangent
section. Indeed, the derivatives of its equation
$x_1q(\zz)+c(\zz)$ are $q(\zz)$ and
$x_1\nabla q(\zz)+\nabla c(\zz)$, and both vanish there.
This proves that the gradients are independent at every common zero.
Such zeros exist in $\PP^{n-3}$ since $n\geq5$.
The two conditions in the statement are open and defined over $K$,
which gives the required $X^\circ$.
\end{proof}

\subsection{Smooth cubic sections of a quadric}

The next lemma supplies the crucial nonsingularity needed to apply the
ten variable theorem. 

\begin{lemma}\label{lem:restriction}
Let $F$ be a field of characteristic zero, let $d\geq2$, and let
$q,c$ be forms over $F$ in $m\geq 2d$ variables, of degrees two and
three, respectively. Suppose that $q$ is nonsingular and
\[
 Y=\{q=c=0\}\subset\PP^{m-1}
\]
is a smooth complete intersection. If $q$ vanishes on a
$d$-dimensional $F$-vector space, then such a space $V$ can be
chosen so that $c|_V$ is a nonsingular cubic form in $d$ variables.
\end{lemma}

\begin{proof}
Choose a $d$-dimensional $F$-space $V_0$ on which $q$ vanishes.
Put $Q=\{q=0\}$ and $G=\operatorname{SO}(q)$. Over $\overline F$,
$G$ is connected and acts transitively on $Q$. Both $\PP(V_0)$ and
$Y$ are smooth and integral over $\overline F$. 
By Kleiman's theorem \cite[Corollary~4(ii)]{kleiman},
for general $g\in G$ the intersection $\PP(gV_0)\cap Y$ is smooth
of dimension
\[
 (d-1)+(m-3)-(m-2)=d-2.
\]
It is nonempty since a cubic on $\PP^{d-1}$ has a geometric zero
when $d\geq2$. Equivalently, the cubic $\xx\mapsto c(g\xx)$ on
$V_0$ has nonzero discriminant for general $g$. In an $F$-basis of
$V_0$ this discriminant is a regular function of $g$ defined over
$F$. Its nonvanishing therefore gives a nonempty open subset of $G$
defined over $F$. 
Since $F$ has characteristic zero, the set $G(F)$ is
Zariski dense in $G$, as explained in 
\cite[Corollary~13.3.10]{springer}.
Thus we can choose $g\in G(F)$ in
this open subset and take $V=gV_0$.
\end{proof}

\subsection{Conclusion}

We  record the results
that we shall need about the arithmetic of  quadratic and  cubic forms.

\begin{proposition}\label{prop:cubic-inputs}
Let $k$ be a number field.
\begin{enumerate}
\item 
For every integer $d\geq1$, 
a quadratic form over $k$ in at least $2d+1$ variables vanishes
on a $d$-dimensional vector space over some quadratic extension $L/k$,
with $L$ totally imaginary.
\item A nonsingular cubic form over $k$ in at least ten variables has
a nontrivial zero over $k$.
\item A smooth cubic hypersurface of dimension at least two with a
$k$-point has Zariski dense $k$-points.
\end{enumerate}
\end{proposition}

\begin{proof}
The first assertion is \cite[Theorem~3]{BDL}, and the second is
\cite[Theorem~1.1]{BV}. The third follows from unirationality
\cite[Theorem~1]{Kollar}.
\end{proof}

\begin{proof}[Proof of Theorem~\ref{thm:cubic}]
Let $X\subset\PP^{n-1}_k$ be a smooth cubic hypersurface over a number field $k$, with
$n\geq23$.
Proposition~\ref{prop:cubic-inputs}(2) gives $X(k)\ne\varnothing$, and
part~(3) allows us to choose $a\in X^\circ(k)$, with $X^\circ$
as in Lemma~\ref{lem:base-point}. Choose coordinates as in 
\eqref{eq:tangent-section} over $k$. We obtain a nonsingular quadratic
form $q$  and a cubic form $c$,
both in $n-2\geq 21$ variables,
 such that
\[
 Y_a=\{q=c=0\}\subset\PP^{n-3}
\]
is a smooth complete intersection.

By Proposition~\ref{prop:cubic-inputs}(1), with $d=10$, there is a
quadratic extension $L/k$ over which $q$ vanishes on a $10$-dimensional
vector space $V$. Lemma~\ref{lem:restriction}, applied over $L$, shows that
this space can be chosen so that $c|_V$ is nonsingular. Applying
Proposition~\ref{prop:cubic-inputs}(2) over $L$ gives $\mathbf0\ne\zz\in V$ such that $q(\zz)=c(\zz)=0$.
The line
\[
 \ell=\{[s:0:t\zz]:[s:t]\in\PP^1_L\}
\]
is thus  contained in $X_L$, by \eqref{eq:tangent-section}, and
passes through $a=[1:0:\cdots:0]\in X(k)$. 
Since $\ell$ passes through the $k$-point $a$,
Wooley's descent lemma \cite[Lemma~2.2]{trevor}
supplies a $k$-line on $X$.
\end{proof}

\end{document}